\documentclass{amsproc}

\usepackage{amssymb,xypic}
\usepackage[dvips]{graphicx}

\newtheorem{theorem}{Theorem}[section]

\theoremstyle{definition}
\newtheorem{definition}[theorem]{Definition}
\newtheorem{eg}[theorem]{Example}
\theoremstyle{remark}
\newtheorem{remark}[theorem]{Remark}

\numberwithin{equation}{section}

\let\oT=\longleftarrow
\let\To=\longrightarrow
\let\into=\hookrightarrow
\newcommand\Ltimes{\mathbin{\buildrel{\mathbf L}\over\otimes}}
\newcommand\OO{\mathcal O}
\newcommand\A{\mathcal A}
\newcommand\B{\mathcal B}
\newcommand\C{\mathbb C}
\newcommand\R{\mathbf R\hspace{1pt}}
\newcommand\RR{\mathcal R}
\newcommand\F{\mathcal F}
\newcommand\Z{\mathbb Z}
\newcommand\res{\arrowvert_}
\newcommand\udot{^{\scriptscriptstyle\bullet}}
\newcommand\ldot{_{\scriptscriptstyle\bullet}}
\newcommand\ve{^\vee}
\newcommand\Hom{\mathrm{Hom\,}}

\begin{document}

\title{Derived categories for the working mathematician}

\author{R.\,P. Thomas}
\address{Mathematical Institute, 24--29 St Giles', Oxford OX1 3LB. UK}
\email{thomas@maths.ox.ac.uk}
\thanks{The author is supported by Hertford College, Oxford}

\date{6 January 99}

\begin{abstract}
It is becoming increasingly difficult for geometers and even
physicists to avoid papers containing phrases like ``triangulated
category'', not to mention derived functors. I will give some
motivation for such things from algebraic geometry, and show how the
concepts are already familiar from topology. This gives a
natural and simple way to look at cohomology and other scary concepts
in homological algebra like Ext, Tor, hypercohomology and spectral
sequences.
\end{abstract}

\maketitle

\section{Introduction}
I should begin by apologising for the title of this talk
\cite{Ma}; while the title is intended for mathematicians the talk itself
is aimed at least as much at physicists. Kontsevich's mirror symmetry
programme \cite{K} and the mathematical description of D-branes have
brought triangulated categories into mainstream string theory and
geometry (there are now even papers in which the reference [Ha] means [RD]
rather than [AG]). Since they have such a fearsome reputation (probably
mainly due to the references being in French), and since I will need
them for my second talk, as presumably will other speakers, I wanted
to motivate this beautiful piece of homological algebra. This
motivation is of course what guided the creators of derived
categories (principally Verdier, or as he is traditionally
known in this context, Grothendieck's student Verdier) but, equally of
course, is not written down. For the technical details of the theory the
working mathematician should consult the excellent \cite{GM}.

The main idea of derived categories is simple: work with complexes
rather than their (co)homology. We will take simple examples from
algebraic geometry to demonstrate why one might want to do this, then
examples from algebraic topology to show that the ideas and
structure are already familiar. (The link between the two subjects is
sheaves, but we shall not pursue this.) This also gives us a
way of using pictures and techniques from topology in algebraic
geometry. Section 1 begins with chain complexes in algebraic geometry
and topology, and why
they are preferable to (co)homology. Thus we would like to consider the
natural invariant of a space/sheaf/etc. to be a complex rather than
its homology, so to do this we must identify which complexes should be
considered isomorphic -- this leads to the notion of quasi-isomorphism
that is the subject of Section 2. Setting up the general (abstract) theory
is where categories creep in, in Section 3, leading to cones, triangles,
and triangulated categories in Section 4. Finally Section 5 shows how
more classical homological algebra, in particular derived functors,
fit into this framework; this makes them more transparent and leads
to easier proofs and explanations of standard results.

\textbf{Acknowledgements.} I learnt a lot of homological algebra from
Brian Conrad, Mikhail Khovanov and Paul Seidel. I would like to
thank members of the CALF (junior
\textbf{C}ambridge-\textbf{O}xford-\textbf{W}arwick) algebraic
geometry seminar for enduring a trial run of a preliminary version of
this talk, and the organisers of the 1999 Harvard Winter School on
Mirror Symmetry for inviting me to take part.

\section{Chain complexes} \label{cx}
\subsection*{Chain complexes in topology}

Chain complexes are familiar in algebraic topology. A
simplicial complex $X$ (such as a triangulation of a manifold) gives
rise to two chain complexes. Letting $C_i$ be the free abelian group
generated by the $i$-simplices in the space, and $C^i$ its dual,
\begin{equation} \label{co}
C_i=\Z\{i\text{-simplices}\}, \qquad C^i=\Hom(C_i,\Z),
\end{equation}
we have the complexes
$$
C\ldot=\ldots\to C_i\stackrel{\partial_i}{\To}C_{i-1}\to\ldots,
\qquad C\udot=\ldots\to C^{i-1}\stackrel{d^i}{\To}C^i\to\ldots.
$$
Here $\partial$ is the boundary operator taking a simplex to its
boundary, $d$ is its adjoint, and $\partial^2=0=d^2$. Thus we can form
the homology and cohomology groups
$$
H_i(X)=\frac{\mathrm{ker\ }\partial_i}{\mathrm{im\ }\partial_{i+1}},
\qquad H^i(X)=\frac{\mathrm{ker\ }d^{i+1}}{\mathrm{im\ }d^i}.
$$
Simplicial maps $X\to Y$ induce chain maps (that is, they commute with
the (co)boundary operators) on the corresponding chain complexes,
inducing maps on the (co)homology groups.

However, the (co)homology $H_*\ (H^*)$ contains less information than
the complexes $C\ldot,\ C\udot$, for instance Massey products. There
exist spaces with the same homology $H_*$ but different homotopy type,
making $H_*$ a limited invariant of homotopy type. But $C\ldot$ is a
very powerful ``invariant'', at least for simply connected spaces
(to which we shall confine ourselves in this talk), due to the
Whitehead theorem. This states that the underlying topological spaces
$|X|,\,|Y|$ of simplicial complexes $X$ and $Y$ (both simply
connected) are homotopy equivalent if and only if there are maps of
simplicial complexes
$$
\spreaddiagramcolumns{-1.2pc}
\spreaddiagramrows{-2.2pc}
\diagram
\qquad & Z\ddlto\ddrto &&&& C\ldot^Z \ddlto\ddrto \\
&&& \mathrm{\ inducing\ chain\ maps\ } \\
X && Y && C\ldot^X & \quad\quad & C\ldot^Y
\enddiagram 
\vspace{-23pt} $$
\begin{equation} \label{wh} \vspace{25pt} \end{equation}
\emph{inducing isomorphisms on homology} $H_*(X)\stackrel{\,\backsim}{\oT}
H_*(Z)\stackrel{\sim\,}{\To}H_*(Y)$. (The reason for
the appearance of $Z$ is that we may need to subdivide $X$ -- thus
giving a $Z$ with a simplicial approximation $Z\to X$ to the
identity map $|Z|\stackrel{\sim\,}{\To}|X|$ -- before we can get a
simplicial approximation $Z\to Y$ to a given homotopy equivalence
$|X|\to|Y|$.)

Thus $C\ldot$ contains as much information as we could hope for, and
we would like to think of it as an invariant of the topological space
$|X|$. We shall see how to do this later; for now we can
think of it as an invariant of the triangulation (simplicial complex).\\

Another place in topology where the advantage of using complexes
rather than their homology is already familiar is the construction of
the dual, cohomology, theory as above. Namely we \emph{do not} define
$H^*$ to be the dual of $H_*$ in the sense that
$$
H^*\ne\Hom(H_*,\Z).
$$
This is because $\Hom(\ .\ ,\Z)$ destroys torsion information
and its square is not the identity (the double dual in this sense is
not the original $H_*$). Instead we take $\Hom(\ .\ ,\Z)$ at
the level of chain complexes as in (\ref{co}). Then no information is
lost, and the double dual is the original complex.

For instance applying $\Hom(\ .\ ,\Z)$ to the complex
$$
\Z\stackrel{2}{\To}\Z \qquad H_1=0,\ H_0=\Z/2,
$$
gives the dual complex
$$
\Z\stackrel{2}{\oT}\Z \qquad H^1=\Z/2,\ H^0=0.
$$
Thus the torsion $\Z/2$ has not been lost on dualising (as it would
have been using $\Hom(H_*,\Z)$), it has just moved in degree.

So we have seen two examples of what will be the theme of this talk,
namely \\
\begin{center} \renewcommand\arraystretch{1.5} \begin{tabular}{|c|}
\hline Complexes good, (co)homology bad \\ \hline
\end{tabular} \end{center} \vspace{3mm}

\subsection*{Chain complexes in algebraic geometry}

It might seem strange to consider chain complexes in algebraic
geometry at all. The main place they arise is in resolutions of
sheaves, which we
describe now. Coherent sheaves often arise naturally as kernels or
cokernels of maps, i.e. as the cohomology of the 2-term complex made
from the map. For instance if $D\subset X$ is a divisor in a smooth
algebraic variety $X$, it corresponds to a line bundle $L$, and a
section $s\in H^0(L)$ vanishing on $D$ (we shall confuse $L$ and its
sheaf of sections $L$). This gives us the standard exact sequence
\begin{equation} \label{div}
0\to L^{-1}\stackrel{s\,}{\to}\OO_X\to\OO_D\to0,
\end{equation}
where $\OO_X$ is the structure sheaf of $X$ (sections of the trivial
line bundle on $X$) and $\OO_D$ is the structure sheaf of $D$ pushed
forward to $X$ (extending by zero -- this is a torsion sheaf
concentrated on $D$). Thus $\OO_D$ is the cohomology of the complex
$\{L^{-1}\stackrel{s\,}{\to}\OO_X\}$. (Notice this is the cohomology
of a complex of sheaves, and as such is a sheaf, and should not be
confused with sheaf cohomology, which is vector-space-valued.)

Similarly if $Z\subset X$ is a codimension $r$ subvariety, the transverse
zero locus of a regular section $s\in H^0(E)$ of a rank $r$ vector
bundle $E$, the exact sequence (Koszul complex)
$$
0\to\Lambda^rE^*\to\Lambda^{r-1}E^*\to\ldots\to E^*\to\OO_X,
$$
where each arrow is given by interior product with $s$, has cokernel
$\OO_Z$ by inspection. Thus $\OO_Z$ is the cohomology of this complex.

We have actually gained something here: we have replaced nasty
torsion sheaves $\OO_D,\ \OO_Z$ by nicer, locally free sheaves
(i.e. vector bundles) on $X$. In general one can consider such
\emph{resolutions}, replacing arbitrary sheaves $\F$ by complexes
$F\udot$ of sheaves that are ``nicer'' in some way,
$$
F\udot\to\F \qquad\mathrm{or}\qquad \F\to F\udot, 
$$
and now use the complex $F\udot$ instead of its less manageable
cohomology $\F$. We think of the nasty looking curly $\F$ being made
up from the nicer straight $F$\,s as $F^1-F^2+F^3-\ldots$, where the
sense in which we subtract the $F^i$s is given by the maps between
them in the resolution, and the $F^i$ are the building blocks: the
generators of $\F$ form $F^1$, the relations $F^2$, relations amongst
the relations $F^3$, and so on.

Of course the resolution may not arise naturally in general
and we must pick one; this is directly analogous to picking a
(non-canonical) triangulation of a topological space as we shall see
in Section \ref{qi}, and is the problem that derived categories resolve.
The sense in which the sheaves $F\udot$ are ``nicer'' will be tackled
in general in Section \ref{df}, but the following is a simple example.

\subsection*{Intersection theory via sheaves}
We concentrate on the easiest case of intersecting two divisors
$D_1,\ D_2$ in a smooth complex surface $X$. These correspond to line
bundles $L_i$ and sections $s_i\in H^0(L_i)$ with zero locus
$D_i$. The $D_i$ have homology classes in $H_2(X)$ that we may
intersect, or dually we can consider $c_1(L_i)\in H^2(X)$. In terms of
sheaf theory the first corresponds to tensoring structure sheaves: if
$D_1,\ D_2$ intersect transversely then we have $\OO_{D_1\cap
D_2}=\OO_{D_1}\otimes\OO_{D_2}$. The second corresponds to using the
resolution (\ref{div}) for $\OO_{D_1}$ and tensoring \emph{that} with
$\OO_{D_2}$, as the following table shows: \\
\renewcommand\arraystretch{2}
$$
\begin{array}{|c||c|l|}
\hline \mathrm{Data} & \mathrm{Divisor}\,\ D_i & s_i\in H^0(L_i)
\mathrm{\ such\ that\ }s_i^{-1}(0)=D_i \\ \hline
\mathrm{Alg\ geom} & \OO_{D_i} & \{L_i^{-1}\stackrel{s_i\,}{\To}
\OO_X\} \\ \hline
\mathrm{Topology} & [D_i]\in H_2(X) & c_1(L_i)\in H^2(X) \\ \hline
\mathrm{Intersection} & [D_1]\,.\,[D_2] & c_1(L_1)\cup c_1(L_2)=
\langle c_1(L_1),D_2\rangle=c_1(L_1\res{D_2}) \\ \hline
& \OO_{D_1}\otimes\OO_{D_2} & c_1(L_1\res{D_2})=(s_1\res{D_2})^{-1}(0),
\mathrm{\ i.e.\ restrict\ } L_1 \vspace{-4mm} \\
\mathrm{Transverse\ case} & =\OO_{D_1\cap D_2} & \mathrm{to\ }
D_2\mathrm{\ and\ take\ zeros\ of\ its\ section\ }s_1\res{D_2},
\vspace{-4mm} \\
&& \mathrm{i.e.\ take\ cokernel\ of\ }\{L_1^{-1}\res{D_2}\stackrel
{s_1\res{D_2}}{\To}\OO_{D_2}\} \\ \hline
\mathrm{Non\ transverse} & \OO_D\otimes\OO_D &
\{L^{-1}\res D\stackrel{s\res D=0\,}{\To}\OO_D\}.\mathrm{\ So\ we\ 
still\ see} \vspace{-4mm} \\
\mathrm{case;\ e.g.} & =\OO_D & L\ \mathrm{on}\
D,\mathrm{\ just\ with\ the\ section\,}=0. \vspace{-4mm} \\
D_1=D_2=D && \\ \hline
\end{array}
$$ \\

So restricting $L_1$ to $D_2$ corresponds to tensoring
$\{L_1^{-1}\to\OO_X\}$ with $\OO_{D_2}$, and in the transverse case
the cokernel of this is just $\OO_{D_1\cap D_2}$. Thus the way to
pass from the right hand column to the left is to take cokernels;
this corresponds to the fact that tensoring with $\OO_{D_2}$ is
\emph{right exact}: tensoring the exact sequence (\ref{div}) with
$\OO_{D_2}$ gives a sequence in which the final three arrows are still
exact.

What the table shows is that while tensoring $\OO_{D_1}$ with
$\OO_{D_2}$ gives the correct answer for $D_1$ and $D_2$ transverse,
in the non-transverse case it does not ($\OO_D$ is the structure sheaf
of $D\cap D$, but not of the correct topological intersection of $D$
with itself). However tensoring instead $\{L_1^{-1}\to\OO_X\}$ with
$\OO_{D_2}$ always gives the right answer: taking the cokernel
of the resulting complex gives the same $\OO_{D_1}\otimes\OO_{D_2}$ as
before, but in the non-transverse case we get more -- we still have
the line bundle $L_1\res{D_2}$ and so the intersection information (we
need only take its first Chern class); by changing the section
$L_1^{-1}\res{D_2}\stackrel{s_1\res{D_2}}{\To}\OO_{D_2}$ from
one which may be identically zero on (components of) $D_2$ to a
transverse section we get the correct intersection. More generally we
take the divisor of $L_1\res{D_2}$ on $D_2$. If $D_1$ can be moved to be
transverse to $D_2$ this gives the same answer; if $L_1\res{D_2}$ has
no sections this cannot be achieved but our method still gives the correct
intersection product, \emph{inside} the scheme-theoretic
intersection $D_1\cap D_2$ as well.

The moral is that
$$
\OO_D\otimes\OO_D \qquad\text{``should\ be''}\qquad
\{L^{-1}\stackrel{s\,}{\to}\OO_X\}\otimes\OO_D=\{L^{-1}\res
D\stackrel{0\,}{\to}\OO_D\}.
$$
This still has cokernel $\OO_D\otimes\OO_D=\OO_D$, but now (because
the intersection was not transverse) has kernel too, namely
$L^{-1}\res D$, containing all the intersection information. The
difference of the first Chern classes of these sheaves on $D$ is
precisely the self intersection $D\,.\,D$.

So just as the dual of $H_*$, in the case of simplicial complexes in
Section \ref{cx}, ``should be'' given by applying $\Hom(\ .\ ,\Z)$
not to $H_*$ but to the chain
complex, here we apply $\otimes\OO_D$ to the complex of locally free
(and so better behaved) sheaves rather than $\OO_D$. This is the
prototype of a derived functor which will be dealt with systematically
in Section \ref{df}. The kernel $L^{-1}\res D$ above will be the first
derived functor Tor$_1(\OO_D,\OO_D)$ of $\otimes\OO_D$.

Thus we see an example where it is beneficial to consider complexes of
sheaves rather than just single (cohomology) sheaves. Having applied
$\otimes\OO_{D_2}$ we get a genuine complex, potentially with
cohomology in more than one degree, i.e. it is not simply the
resolution of a single sheaf. But again we should not now pass
to its cohomology, as we may want to intersect with further cycles, by
tensoring the complex with some $\OO_{D_3}$ if
$X$ is higher dimensional, for instance. Complexes good, cohomology
bad, after all.

So having motivated considering \emph{all} complexes (rather than just
resolutions of sheaves) we will now set
about working with them, forgetting all about
derived functors until Section \ref{df}. Replacing sheaves by
complexes of which they are the cohomology, i.e. by resolutions, we
come across the problem mentioned earlier: how do we
pick a resolution functorially\,? Again there is an analogous issue in
topology.

\section{Quasi-isomorphisms} \label{qi}
\subsection*{Quasi-isomorphisms in topology}

We saw in Section \ref{cx} that a good invariant of a (simply
connected) topological space $|X|$ underlying a simplicial complex $X$
is the simplicial chain complex $C\ldot^X$ which determines the
homotopy type of $|X|$ by the Whitehead theorem. The problem is
functoriality: how to pick a triangulation of $|X|$ canonically, to
pass from the topological space to a complex. The standard
mathematical trick is to consider all at once on an equal footing, for
instance by making them all isomorphic. As described in Section
\ref{cx}, different
triangulations of a space (yielding different chain complexes
$C\ldot,\,D\ldot$) may have no simplicial map between them, but
taking finer subdivisions and using simplicial approximations we can
find a third chain complex $E\ldot$ fitting into the diagram
$$
\spreaddiagramcolumns{-1pc}
\spreaddiagramrows{-1pc}
\diagram
& E\ldot \dlto\drto \\ C\ldot && D\ldot
\enddiagram
$$
where both maps are \emph{quasi-isomorphisms} -- they are chain maps
inducing isomorphisms on homology. We extend this to be an equivalence
relation; thus two chain complexes are quasi-isomorphic if they can be
related by a sequence of quasi-isomorphisms of the above type.
So quasi-isomorphic complexes have the same homology, but the
converse does \emph{not} hold; for instance the complexes
$$
\C[x,y]^{\oplus2}\stackrel{(x,y)\,}{\To}\C[x,y] \qquad\mathrm{and}\qquad
\C[x,y]\stackrel{0\,}{\to}\C
$$
have the same homology but are not quasi-isomorphic. Quasi-isomorphism
is exactly the equivalence relation we want on complexes: by the
Whitehead theorem (\ref{wh}), $|X|$ and $|Y|$ are homotopy equivalent
if and only if there is a $Z$ inducing a quasi-isomorphism
$$
\spreaddiagramcolumns{-1pc}
\spreaddiagramrows{-1pc}
\diagram
& C^Z\ldot \dlto\drto \\ C^X\ldot \ar@{-->}[rr] && C^Y\ldot.
\enddiagram \vspace{-24pt} $$
\begin{equation} \label{dot} \vspace{26pt} \end{equation}
Thus we would like to think of quasi-isomorphic complexes as
isomorphic. Though there may not be a map between them (as there may be no
simplicial map between two different triangulations of the same space)
we pretend there is by putting one in by hand, putting in the
dotted arrow in (\ref{dot}) even if it does not exist as a genuine
map. We consider the complexes to be
isomorphic since the underlying spaces are, after all.

In particular homotopy equivalences are
quasi-isomorphisms. Algebraically this means that if there is an
$s:\,C\ldot\to D_{\scriptscriptstyle\bullet+1}$ such that
$f-g=\partial_D\circ s+s\circ\partial_D$ for two chains maps
$f,\,g:C\ldot\to D\ldot$, then $f$ and $g$ induce the same map on
homology.

\subsection*{Quasi-isomorphisms in algebraic geometry}

The cokernel map (\ref{div})
$$
\{L^{-1}\to\OO_X\}\To\OO_D
$$
is of course a quasi-isomorphism (it induces an isomorphism on
cohomology). Similarly any resolution $F\udot\to\F$ (or $\F\to
F\udot$) gives a quasi-isomorphism $\{F\udot\}\to\F$ (or $\F\to
\{F\udot\}$). As in the topological case we want to consider these as
equalities (more strictly, isomorphisms) $\F\cong\{F\udot\}$, to make
choice of a resolution functorial. So again we need to introduce
inverse arrows $\{F\udot\}\dashleftarrow\F$ (or $\F\dashleftarrow
\{F\udot\}$), which is what we turn to now.

\section{Category theory}

In general we would like to replace \emph{objects} (e.g. sheaves, vector
bundles, abelian groups, modules, etc.) by complexes of objects that
are \emph{``better behaved''} for some particular \emph{operation}
such as $\otimes A,\ \Hom(\ .\ ,A),\ \Hom(A,\ .\ ),\ \Gamma(\ .\
),\,\ldots$ (better behaved in the sense that \emph{no information is
lost} when the operation is applied to such objects). Sometimes there
is a canonical such complex, in general we want to make its choice
\emph{natural}.

The words in italics are meant to be suggestive of category theory;
abelian categories are the natural setting for the general theory, with
the above words being replaced by, respectively, objects of a
category, acyclic objects, functor, exact, functorial. We shall now
define such things; those scared of the word category should think of
the (category of) sheaves on a complex manifold.

\begin{definition}
An additive category is a category $\A$ such that
\begin{itemize}
\item Each set of morphisms $\Hom(A,B)$ forms an abelian group.
\item Composition of morphisms distributes over the addition of
  morphisms given by the abelian group structure,
  i.e. $f\circ(g+h)=f\circ g+f\circ h$ and $(f+g)\circ h=f\circ
  h+g\circ h$.
\item There exist products (direct sums) $A\times B$ of any two
  objects $A,\,B$ satisfying the usual universal properties (see
  e.g. \cite{Ma}).
\item There exists a zero object $0$ such that $\Hom(0,0)$ is the zero
  group (i.e. just the identity morphism). Thus
  $\Hom(0,A)=0=\Hom(A,0)$ for all $A$, and the unique zero morphism
  between any two objects is the one that factors through the zero object.
\end{itemize}
An abelian category is an additive category that also satisfies
\begin{itemize}
\item All morphisms have kernels and cokernels such that monics
  (morphisms with zero kernel) are the kernels of their cokernel maps
  and epis (zero cokernel) are the cokernels of their kernels.
\end{itemize}
\end{definition}

The kernels and cokernels mentioned above are defined in a general
category via obvious universal properties (e.g. a kernel of
$f:\,A\to B$ is a $K\to A$ through which any $g:\,C\to A$
factors uniquely if and only if $f\circ g=0$). We refer to \cite{Ma}
for such things; we shall only be concerned with concrete categories
(those whose objects are sets) and then kernels and cokernels will be
the usual ones, and the final axiom will be immediate.

So in an abelian category we can talk about exact sequences and
\emph{chain complexes}, and cohomology of complexes. Additive functors
between abelian categories are \emph{exact} (respectively left or
right exact) if they preserve exact sequences (respectively short
exact sequences $0\to A\to B\to C$ or $A\to B\to C\to0$).

\begin{definition}
The bounded derived category $D^b(\A)$ of an abelian category $\A$ has
as objects bounded (i.e. finite length) $\A$-chain complexes, and
morphisms given by chain maps with quasi-isomorphisms inverted as
follows (\cite{GM} III 2.2). We introduce morphisms 
$f$ for every chain map between complexes $f:\,X_f\to
Y_f$, and $g^{-1}:\,Y_g\to X_g$ for every quasi-isomorphism
$g:\,X_g\stackrel{\sim\,}{\to}Y_g$. Then form all products of these
morphisms such that
the range of one is the domain of the next. Finally identify any
combination $f_1f_2$ with the composition $f_1\circ f_2$, and
$gg^{-1}$ and $g^{-1}g$ with the relevant identity maps id$_{Y_g}$ and
id$_{X_g}$.
\end{definition}

\begin{remark}
Similarly one can define the unbounded derived category, and the
categories $D^+(\A),\ D^-(\A)$ of bounded below and above complexes
respectively. We shall use $D(\A)$ to mean one of the four such
derived categories.
\end{remark}

Morphisms in $D(\A)$ are represented by roofs
$$
\diagram
& Z\udot \dlto_s\drto^g \\ X\udot \ar@{-->}[rr]^f && Y\udot,
\enddiagram
$$
where $s$ is a quasi-isomorphism, and we set $f=gs^{-1}$ (it is clear
any morphism is a composition of such roofs; that one will suffice is
proved in \cite{GM} III 2.8). This
``localisation'' procedure of inverting quasi-isomorphisms has some
remarkable properties (for instance we shall see that
homotopic maps are identified with each other to give
the same morphism in the derived category). Although they give $D(\A)$
the structure of an additive category, we will see it does not have
kernels or cokernels. As ever we go back to the topology to see why
not, what they are replaced by, and what the structure of $D(\A)$ is
(since it is not that of an abelian category).

\section{Cones and triangles}

When working with topological spaces (or simplicial or cell complexes)
up to homotopy there is no notion of kernel or cokernel. In fact the
standard cylinder construction shows that any map
$f:\,X\to Y$ is homotopic to an inclusion
$X\to\,$cyl$\,(f)=Y\sqcup(X\times[0,1])/f(x)\sim(x,1)$, while the path
space construction shows it is also homotopic to a fibration.

For some \emph{fixed} maps $f:\,X\to Y$, rather then equivalence classes
of homotopic maps, we can make sense of the kernel
(the fibre of $f$ if it is a fibration) or
cokernel ($Y/X$, the space with the image of $X$ collapsed to a point,
if $f$ is a cofibration, which means an inclusion for our
purposes). In general of course neither makes sense, but there is
something which acts as both, namely cones and the Dold-Puppe
construction.

{The cone $C_f$ on a map $f:\,X\to Y$ is the space formed from
$Y\sqcup(X\times[0,1])$ by identifying $X\times\{1\}$ with its
image $f(X)\subset Y$, and collapsing $X\times\{0\}$ to a point. It
fits into the sequence of maps}
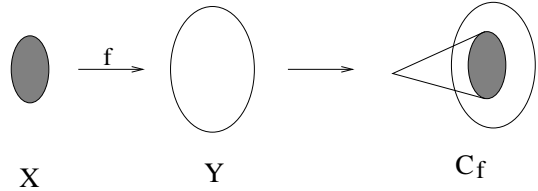
\begin{figure}[h]
\begin{center}
\includegraphics[angle=0,width=7cm]{cone.eps}
\end{center}
\end{figure}

{It is clear that this can act as a cokernel, in that if $X\stackrel{f}
{\into}Y$ is an inclusion, then $C_f$ above is clearly homotopy
equivalent to $Y/X$. In fact we can now iterate this process, forming
the cone on the natural inclusion $i:\,Y\to C_f$ to give the sequence:}
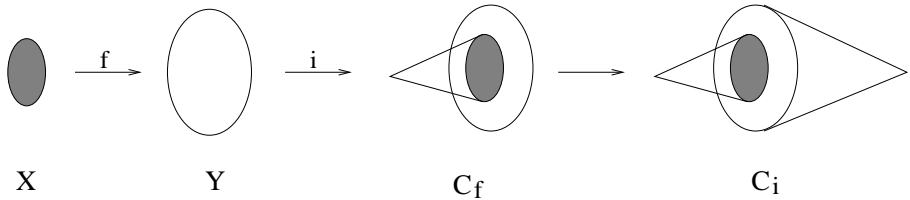
\begin{figure}[h]
\begin{center}
\includegraphics[angle=0,width=12cm]{doldpup.eps}
\end{center}
\end{figure}

By retracting $Y$ to the right hand apex of $C_i$ in the last term
of the above diagram, we see that $C_i$ is homotopic to $\Sigma X$,
the suspension of $X$. Thus, up to homotopy, we get a sequence
\begin{equation} \label{sp}
X \To Y \To Y/X \To \Sigma X \To \ldots
\end{equation}

Taking the $i$th cohomology $H_i$ of each term, and using the
suspension isomorphism $H_i(\Sigma X)\cong H_{i-1}(X)$ gives a sequence
\begin{equation} \label{les}
H_i(X)\to H_i(Y)\to H_i(Y,X)\to H_{i-1}(X)\to H_{i-1}(Y)\to\ldots
\end{equation}
which is just the long exact sequence associated to the pair
$X\subset Y$ (it is an interesting exercise to check that the map
$H_i(Y/X)\to H_i(\Sigma X)$ induced above is indeed the boundary map
$H_i(Y,X)\stackrel{\partial_*}{\To}H_{i-1}(X)$).

Up to homotopy we can make (\ref{sp}) into a sequence of
simplicial maps, so that taking the associated chain complexes we get
a lifting of the long exact sequence of homology (\ref{les}) to the
level of complexes. Of course it exists for all maps $f$, not just
inclusions, with $Y/X$ replaced by $C_f$.

For instance if $f$ is a fibration, $C_f$ acts as the ``kernel'' or
fibre of the map. The most extreme case is $f:\,X\to$\,point, which
gives $C_f=\Sigma X$, the suspension of the fibre $X$, which to
homology is $X$ shifted in degree by 1.

So $C_f$ acts as a combination of both cokernel and kernel, and as
such gives an easy proof of the Whitehead theorem we have been
quoting: i.e. if $f:\,X\to Y$ is a map inducing an isomorphism of
homology groups of simply connected spaces then the sequence
$$
H_i(X)\to H_i(Y)\to H_i(C_f)\to H_{i-1}(X)\to H_{i-1}(Y)\to\ldots
$$
shows that $H_*(C_f)$=0. Thus, by the Hurewicz
theorem, the homotopy groups $\pi_*(C_f)$ are zero too, making $C_f$
homotopy equivalent to a point by the more famous Whitehead
theorem. It is an easy consequence of this that $f$ is a homotopy
equivalence.

The final thing we note from the topology is that if $X$ and $Y$ are
simplicial complexes, with $f:\,X\to Y$ a simplicial map, then the
cone $C_f$ is naturally a simplicial complex with $i$-simplices being
those in $Y$, plus the cones on $(i-1)$-simplices in $X$. A little
thought shows that this makes the corresponding
(cohomology, for convenience) complex
$$
\renewcommand\arraystretch{1}
C_X\udot\,[\,1\,]\oplus C_Y\udot \quad\mathrm{with\ differential}\quad
d_{C_f}=\left(\!\!\!\begin{array}{cc}
d_X\,[\,1\,] & 0 \\ f & d_Y \end{array}\!\!\right),
$$
where $[\,n\,]$ means shift a complex $n$ places left.

Thus we can define the cone $C_f$ on any map of chain complexes
$f:\,A\udot\to B\udot$ in an
abelian category $\A$ by the above formula, replacing $C_X\udot$
by $A\udot$ and $C_Y\udot$ by $B\udot$. If $A\udot=A$ and $B\udot=B$
are chain complexes concentrated in degree zero then $C_f$ is the
complex $\{A\stackrel{f\,}{\to}B\}$. This has zeroth cohomology
$h^0(C_f)=$\,ker$\,f$, and $h^1(C_f)=$\,coker\,$f$, so combines the
two (in different degrees). In general it is just the total complex of
$A\udot\to B\udot$.

There is an obvious map
$i:\,B\udot\to C_f$, and $C_i$ is, as could be guessed from above,
quasi-isomorphic to $A\udot\,[\,1\,]$. So what we get in a derived
category is not kernels or cokernels, but ``exact triangles''
$$
A\udot\to B\udot\to C\udot\to A\udot\,[\,1\,].
$$
(For a proof see (\cite{GM} III 3.5) -- we should really modify $B$,
replacing it by the quasi-isomorphic cyl\,$(f)$, but once
quasi-isomorphisms are inverted it does not matter.)

Thus we have long exact sequences rather than short exact ones; taking
$i$th cohomology $h^i$ of the above gives the standard long exact sequence
$$
h^i(A\udot)\to h^i(B\udot)\to h^i(C\udot)\to h^{i+1}(A\udot)\to\ldots
$$

Thus $D(A)$ is not an abelian category, it is an example of a
\emph{triangulated category}. This is an additive category with a
functor $T$ (often denoted $[\,1\,]$) and a set of \emph{distinguished
triangles} satisfying a list of axioms. We refer to (\cite{GM} IV 1.1)
for the precise details, but the triangles include, for all objects
$X$ of the category,
$$
X\stackrel{\mathrm{id}\,}{\To}X\to0\to X\,[\,1\,],
$$
and any morphism $f:\,X\to Y$ can be completed to a distinguished
triangle
$$
X\to Y\to C\to X\,[\,1\,].
$$
There is also a derived analogue of the 5-lemma, and a compatibility of
triangles known as the octahedral lemma, which is pretty un\TeX able,
as you might imagine.

Why are we abstracting again\,? It is because this structure that is
present in $D(\A)$, for an abelian category $\A$, is also present in
topology (as we have been hinting all along) without any
underlying abelian category. I.e. the topology described throughout
this talk might lead one to suspect that some category of nice
topological spaces should have just this structure,
and we have certainly not constructed it as a derived category of any
abelian category.

In fact we are not quite there yet, the reason being that the
translation functor of suspension $T=\Sigma$ is not invertible on
spaces. If it did exist it is clear it would be the loop space functor
$\Omega$, as $\Omega$ and $\Sigma$ are adjoints by a standard
argument\footnote{If we denote homotopy classes of maps from $X$ to $Y$
 by $[X,Y]$, then applying the invertible $\Sigma$ to $[X,\Sigma^{-1}Y]$
 gives an isomorphism $[X,\Sigma^{-1}Y]\cong[\Sigma
 X,Y]\cong[X,\Omega Y]$. Choosing $X=\Sigma^{-1}Y$ and $X=\Omega Y$
 gives canonical maps $\Sigma^{-1}Y\leftrightarrows\Omega Y$ that
 induce isomorphisms on homotopy groups (as is seen by choosing
 $X=S^n$). Thus $\Sigma^{-1}Y$ is, up to weak homotopy equivalence,
 $\Omega Y$.}. In particular $\Sigma X\simeq\Omega^{-1}X$ and we would
be able to deloop spaces, making them infinite loop spaces.

So it is not surprising to find that \emph{spectra} (essentially
infinite loop spaces, plus all their loopings and deloopings to keep
track of homotopy in negative degrees) form a triangulated
category in the way we have been describing (though there are many
technicalities we have bypassed, due to the spaces involved not being
simplicial complexes).

The final thing to note from the topology is that homotopic maps $X\to
Y$ get identified when we invert quasi-isomorphisms. By composing with
the homotopy $X\times[0,1]\to Y$ it is clear it is enough to show the
two inclusions $\iota_0,\,\iota_1$ of $X$ into $X\times[0,1]$ (as
$X\times\{0\}$ and $X\times\{1\}$ respectively) are identified. But
they are both right inverses for the projection $p:\,X\times[0,1]\to X$:
$p\iota_i=$\ id. As they are also isomorphisms on homology, when we invert
quasi-isomorphisms they both become identified with $p^{-1}$.

Algebraically, in $D(\A)$, we can mimic the (dual, cohomology) proof
by defining the cylinder of a complex in the way dictated by the
obvious product cell complex structure on the cylinder of a simplicial
complex. Namely consider, for any chain complex $A\udot$,
$$
\renewcommand\arraystretch{1}
\mathrm{cyl}\,(A\udot)=A\udot\oplus A\udot\,[\,1\,]\oplus A\udot,
\quad\mathrm{with\ differential}\quad d=\left(\!\!\!\begin{array}{ccc}
d_A & 0 & 0 \\ -\mathrm{id} & d_A\,[\,1\,] & \mathrm{id} \\ 0 & 0 & d_A
\end{array}\!\!\!\right).
$$
Set $\iota_1^*,\,\iota_2^*$ to be the projections to the first and
third factors respectively, with $p^*$ the sum of the corresponding
two inclusions. These are all chain maps, and in fact
quasi-isomorphisms, with $\iota_i^*p^*=\ $id. Thus
$\iota_1^*,\,\iota_2^*$ are identified with $(p^*)^{-1}$ in $D(\A)$
and so also with each other.

\section{Derived functors} \label{df}

Now we can go back to (derived) functors and replace arbitrary
complexes by resolutions of objects which are suited to the functor
concerned. We deal with left exact functors; right exact functors are
similar.

\begin{definition}
Let $\A$ and $\B$ be abelian categories. A class of objects
$\RR\subset\A$ is \emph{adapted} (\cite{GM} III 6.3) to a left exact functor
$F:\A\to\B$ if
\begin{itemize}
\item $\RR$ is stable under direct sums,
\item $F$ applied to an acyclic complex in $\RR$ (that is, a complex
  with vanishing cohomology) is acyclic, and
\item any $A\in\A$ injects $0\to A\to R$ into some $R\in\RR$.
\end{itemize}
Let $K^+(\RR)$ be the category of bounded below chain complexes in $\RR$
with morphisms homotopy equivalence classes of chain maps. Then
\cite{GM} inverting quasi-isomorphisms in $K^+(\RR)$ gives a category
equivalent to $D^+(\A)$.
\end{definition}

The final statement is a fancy way of saying we have finally reached
our goal -- we can functorially replace (i.e. resolve) any
$\A$-complex by a quasi-isomorphic $\RR$-complex using the conditions
of the definition. Examples of such $\RR\subset\A$ include $\{$injective
sheaves$\}\subset\{$quasi-coherent sheaves$\}$, or for right exact
functors, $\{$projectives modules over a ring$\}\subset\{$all
modules$\},\ \{$locally free sheaves over an affine variety or
scheme$\}\subset\{$coherent sheaves$\},\ \{$flat
sheaves$\}\subset\{$quasi-coherent sheaves$\}$ if $F$ is tensoring
with a sheaf, etc.

So instead of applying $F$ to arbitrary complexes, we apply it only to
$\RR$-complexes.

Thus we define the \emph{right derived functor} $\R F$ of $F$ to be
the composition $D^+(\A)\to K^+(\RR)/\mathrm{q.\,i.\,}\stackrel{F\,}
{\To}D^+(\B)$. (The more classical right derived functors $R^iF$ are
the cohomology $R^iF=h^i(\R F)$, but taking cohomology is bad, we now
know.) We are glossing over some details here; really $\R F$ should be
defined by some universal property (\cite{GM} III 6.11) to make it
independent of $\RR$ up to canonical isomorphism. Similarly we can
define the left derived functor $\mathbf LF$ of a right exact $F$.

This gives an \emph{exact functor} $\R F:\,D^+(\A)\to D^+(\B)$,
i.e. it takes exact triangles to exact triangles. In particular taking
cohomology gives the classical long exact sequence of derived functors
of the form
$$
\ldots\to R^iF(A)\to R^iF(B)\to R^iF(C)\to R^{i+1}F(A)\to\ldots
$$

For instance the right derived functor of the left exact global
sections functor $\Gamma:\,\{$\,Sheaves\,$\}\to\{$\,Vector
spaces\,$\}$ is just sheaf cohomology $\R\Gamma=\R H\udot$, with
cohomology (of the complex of vector spaces) the standard sheaf
cohomology $H^i$. Then the above sequence becomes the long exact
sequence in cohomology.

There are two main advantages of this approach. Firstly that we have
managed to make the \emph{complex} $\R F(A)$, rather than its less
powerful cohomology $R^iF(A)$, into an invariant of $A$, unique up to
quasi-isomorphism. Secondly, the derived functor has simply
become the original functor applied to complexes (though not
arbitrary ones, they have to be in $\RR$). This gives easier and more
conceptual proofs for results about derived functors that usually require
complicated double complex, spectral-sequence type arguments. We give
some examples in sheaf theory, skating over a few technical conditions
(issues about boundedness of complexes and resolutions that are
certainly not a problem on a smooth projective variety; for precise
statements see \cite{HRD}):

\begin{eg}
Tensor product of sheaves is symmetric, thus its derived
functor $\Ltimes$, and its homology Tor$_i$, are symmetric:
Tor$_i(A,B)\cong$\,Tor$_i(B,A)$. Here we simply tensor complexes of
flat sheaves (for instance, locally free sheaves).
\end{eg}

\begin{eg}
$\R\mathcal Hom$ of sheaves, being just local $\mathcal
Hom$ on complexes, can be defined by resolving either the first
variable by locally frees or the second by injectives.
\end{eg}

\begin{eg}
Under some mild conditions, $\R(F\circ
G)\cong\R F\circ\R G$. If we take cohomology before applying $\R F$ we
get an approximation to $\R(F\circ G)$, and the Grothendieck spectral
sequence $R^iF(R^jG)\Rightarrow R^{i+j}(F\circ G)$.

For instance for a morphism $p:\,X\to Y$ the equality $\Gamma_X=
\Gamma_Y\circ p_*$ yields $\R\Gamma_X\cong\R\Gamma_Y\circ\R p_*$, and
so the Leray spectral sequence $H^{i+j}_X\Leftarrow H^i_Y(R^jp_*)$.

Similarly $\Hom=\Gamma\circ\mathcal Hom$ yields $\R\Hom\cong\R
H\udot(\R\mathcal Hom)$ and the local-to-global spectral sequence
Ext$^{i+j}\Leftarrow H^i(\mathcal Ext^j)$.

Hypercohomology $\mathbb H$ is given by the derived functor
of global sections of complexes (so it is nothing but normal sheaf
cohomology for us). Applying the above to $F=\Gamma$ and $G=$\,id we
get the hypercohomology spectral sequence $H^i(h^j(A\udot))\Rightarrow
\mathbb H^{i+j}(A\udot)$.
\end{eg}

\begin{eg}
For $A$ a sheaf or complex of sheaves, denote by $A\ve$ the dual
complex $\R\mathcal Hom(A,\OO)$. Then $\R\mathcal Hom(A,B)\cong B\Ltimes
A\ve$.
\end{eg}

\begin{eg}
We now consider an example on a real $n$-manifold, namely the DeRham
theorem. Let $A^i(\mathbb R)$ denote the sheaf of $C^\infty$ $i$-forms
on a manifold $X$, with $d$ the exterior derivative. Then the
Poincar\'e lemma gives a quasi-isomorphism (resolution of the constant
sheaf $\mathbb R$)
$$
\mathbb R\to\{A^0(\mathbb R)\stackrel{d\,}{\to}A^1(\mathbb R)\stackrel
{d\,}{\to}\ldots\stackrel{d\,}{\to}A^n(\mathbb R)\}.
$$
Because of the existence of partitions of unity the $A^i(\mathbb R)$
sheaves are acyclic for the global sections functor $\Gamma$,
i.e. they have no higher sheaf cohomology (they are what is known as
\emph{fine} sheaves). Thus it is easy to see that we
may apply $\Gamma$ to the resolution to obtain a complex isomorphic to
the derived functor $\R\Gamma(\mathbb R)$ of the constant sheaf
$\mathbb R$, i.e. its sheaf cohomology. This yields
$$
\Omega^0(\mathbb R)\stackrel{d\,}{\to}\Omega^1(\mathbb R)\stackrel
{d\,}{\to}\ldots\stackrel{d\,}{\to}\Omega^n(\mathbb R),
$$
which means that the DeRham complex computes the real cohomology
$H^*(X;\mathbb R)$ of the manifold $X$.

Similarly taking cohomology of the Dolbeault complex or resolution
$$
\C\to\{\OO\stackrel{\partial\,}{\to}\Omega^{1,0}\stackrel
{\partial\,}{\to}\Omega^{2,0}\stackrel{\partial\,}{\to}\ldots
\stackrel{\partial\,}{\to}\Omega^{n,0}\},
$$
on a complex manifold $X$, gives a spectral sequence relating $H^*(X;\C)$
to the Hodge groups $\bigoplus_{i,j}H^i(\Omega^{j,0})=
\bigoplus_{i,j}H^{i,j}(X)$
(the spectral sequence famously degenerates for a K\"ahler manifold).
\end{eg}

Hopefully this talk has shown derived and triangulated categories to
be natural objects, just categories of complexes with quasi-isomorphisms
made into isomorphisms, which we can think of via topological
pictures (modulo some technical details). Unfortunately this is
a bit of a disservice to anyone who comes to work in some $D(\A)$,
where it is important to chase all quasi-isomorphisms,
compatibilities, etc. The main problem is that while cones are defined
up to isomorphism in triangulated categories, choosing them
functorially is not possible. This leads many to think that there
should be some more refined concept still to be worked out.

\bibliographystyle{amsalpha}

\end{document}